\documentclass[journal]{IEEEtran}
\ifCLASSINFOpdf
\else
\fi
\usepackage[pass]{geometry}
\usepackage{fancyhdr,todo}
\usepackage{amsmath}
\usepackage{amsfonts,amssymb,amsthm,cite,float,graphicx}
\usepackage[numbers]{natbib}    

\usepackage{subfigure}
\usepackage[usenames]{color}
\usepackage{pifont}
\usepackage{enumitem}

\usepackage{multirow}
\usepackage{hhline}

\usepackage{listings}
\usepackage{color} 
\definecolor{mygreen}{RGB}{28,172,0} 
\definecolor{mylilas}{RGB}{170,55,241}

\usepackage{algpseudocode,algorithm,algorithmicx}
\newcommand*\Let[2]{\State #1 $\gets$ #2}

\theoremstyle{plain}
\newtheorem{theorem}{Theorem}[section] 

\theoremstyle{definition}
\newtheorem{assumption}[theorem]{Assumption}
\newtheorem{lemma}[theorem]{Lemma}

\usepackage{epstopdf} 

\begin{document}
	
\title{Augmented Lagrangian Optimization \\under Fixed-Point Arithmetic}

\author{Yan~Zhang,~\IEEEmembership{Student Member,~IEEE,}
        Michael~M.~Zavlanos,~\IEEEmembership{Member,~IEEE}
\thanks{The authors are with the Department
of Mechanical Engineering and Material Science, Duke University, Durham, NC, 27708, USA e-mail: yz227@duke.edu; michael.zavlanos@duke.edu.}}


\maketitle

\begin{abstract}
In this paper, we propose an inexact Augmented Lagrangian Method (ALM) for the optimization of convex and nonsmooth objective functions subject to linear equality constraints and box constraints where errors are due to fixed-point data. To prevent data overflow we also introduce a projection operation in the multiplier update. We analyze theoretically the proposed algorithm and provide convergence rate results and bounds on the accuracy of the optimal solution. Since iterative methods are often needed to solve the primal subproblem in ALM, we also propose an early stopping criterion that is simple to implement on embedded platforms, can be used for problems that are not strongly convex, and guarantees the precision of the primal update. To the best of our knowledge, this is the first fixed-point ALM that can handle non-smooth problems, data overflow, and can efficiently and systematically utilize iterative solvers in the primal update. Numerical simulation studies on a utility maximization problem are presented that illustrate the proposed method.
\end{abstract}


\IEEEpeerreviewmaketitle

\section{Introduction}
\label{sec:intro}
Embedded computers, such as FPGAs (Field-Programm-able Gate Arrays), are typically low-cost, low-power, perform fast computations, and for these reasons they have long been used for the control of systems with fast dynamics and power limitations, e.g., automotive, aerospace, medical, and robotics. While embedded computers have been primarily used for low-level control, they have recently also been suggested to obtain real-time solutions to more complex optimization problems \cite{jerez2014embedded,richter2017resource}.
The main challenges in implementing advanced optimization algorithms on resource-limited embedded devices are providing complexity certifications and designing the data precision to control the solution accuracy and the dynamic range to avoid overflow.

While recent methods, such as \cite{jerez2014embedded,richter2017resource} , address these challenges and provide theoretical guarantees, they do so for quadratic problems. In this paper, we provide such guarantees for general convex problems. Specifically, we propose an Augmented Lagrangian Method (ALM) to solve problems of the form
\begin{equation}
	\label{pb:problem}
	\begin{split}
		\min \quad f(x) \;\;\; \text{s.t.} \;\;\; Ax = b, \;\;\; x \in \mathcal{X}
	\end{split}
\end{equation}
on embedded platforms under fixed-point arithmetic, where $x \in \mathbb{R}^n$, $A \in \mathbb{R}^{p \times n}$, $b \in \mathbb{R}^p$ and $f(x)$ is a scalar-valued function. 
The set $\mathcal{X}$ is convex.
ALM falls in the class of first order methods, which have been demonstrated to be efficient in solving problem~\eqref{pb:problem} on embedded computers due to their simple operations and less memory requirements, \cite{richter2017resource}.

Recent work on error analysis of inexact first-order methods is presented in \citet*{devolder2014first,patrinos2015dual,necoara2016iteration}. Specifically, \citet{devolder2014first} proposed a first-order inexact oracle to chracterize the iteration complexity and suboptimality of the primal gradient and fast gradient method. \citet{patrinos2015dual,necoara2016iteration} extend these results to the dual domain.
However, these analyses assume strong convexity of the objective functions. \citet*{nedelcu2014computational,necoara2017complexity} analyze the convergence of ALM using an inexact oracle for general convex problems. The inexactness comes from the approximate solution of the subproblems in ALM. Similar analysis has been conducted in \citet*{rockafellar1976monotone,eckstein2013practical,lan2016iteration}. However, none of above works on ALM consider the error in the multiplier update under fixed-point arithmetic. Moreover, since no projection is used in the multiplier update, the above works cannot provide an upper bound on the multiplier iterates and therefore cannot avoid data overflow. Perhaps the most relevant work to the method proposed here is \citet*{jerez2014embedded}. Specifically, \cite{jerez2014embedded} analyzes the behavior of the Alternating Direction Method of Multipliers (ADMM), a variant of the ALM, on fixed-point platforms. However, the analysis in \cite{jerez2014embedded} can only be applied to quadratic objective functions. Moreover, to prevent data overflow, the proposed method needs to monitor the iteration history of the algorithm to estimate a bound on the multiplier iterates.

Compared to existing literature on inexact ALM, we propose a new inexact ALM that incorporates errors in both the primal and dual updates and contains a projection operation in the dual update. Assuming a uniform upper bound on the norm of the optimal multiplier is known, this projection step can guarantee no data overflow during the whole iteration history. To the best of our knowledge, this is the first work to provide such guaranetees for general convex and non-smooth problems. Furthermore, we show that our proposed algorithm has $O(1/K)$ convergence rate and provide bounds on the achievable primal and dual suboptimality and infeasibility. In general, iterative solvers are needed to solve the subproblems in ALM but the theoretical complexity of such solvers is usually conservative. Therefore, we present a stopping criterion that allows us to terminate early the primal iteration in the ALM while guaranteeing the precision of the primal update. This stopping condition is simple to check on embedded platforms and can be used for problems that are not necessarily strongly convex. 
We note that in this paper we do not provide a theoretical uniform upper bound on the optimal multiplier for general convex problems. Instead, our contribution is to develop a new projected ALM method that relies on such bounds to control data overflow on fixed-point platforms.

The rest of this paper is organized as follows. In Section~\ref{sec:prelim}, we formulate the problem and introduce necessary notations and lemmas needed to prove the main results. In Section~\ref{sec:converg}, we characterize the convergence rate of the algorithm and present bounds on the primal suboptimality and infeasiblity of the solution. In Section~\ref{sec:fp}, we present the stopping criterion for the solution of the primal subproblem under fixed-point arithmetic. In Section~\ref{sec:sim}, we present simulations that verify the theoretical analysis in the previous sections. In Section~\ref{sec:conclude}, we conclude the paper.

\section{Preliminaries}
\label{sec:prelim}
We make the following assumptions on problem~\eqref{pb:problem}.
\begin{assumption}
	\label{assum:prob_form}
	The function $f(x)$ is convex and is not necessarily differentiable. The problem~\eqref{pb:problem} is feasible.
\end{assumption}

\begin{algorithm}[t]
	\caption{Augmented Lagrangian Method}
	\label{alg:exactAL}
	\begin{algorithmic}[1]
		\Require{Initialize $\lambda_0 \in \mathbb{R}^p$, k=0}
		\While{$Ax_{k} \neq b$}
		\State{$x_{k} = \arg \min_{x\in \mathcal{X}} L_\rho(x;\lambda_k)$}
		\State{$\lambda_{k+1} = \lambda_k + \rho(Ax_{k} - b)$}
		\Let{$k$}{$k + 1$}
		\EndWhile
	\end{algorithmic}
\end{algorithm}

The Lagrangian function of problem~\eqref{pb:problem} is defined as
$L(x; \lambda) = f(x) + \langle Ax - b, \lambda \rangle$ and the dual function is defined as
$ \Phi(\lambda) \triangleq \min_{x\in\mathcal{X}} L(x;\lambda)$,  where $\lambda \in \mathbb{R}^n$ is the Lagrangian multiplier \citep*{ruszczynski2006nonlinear} and $\langle \cdot , \cdot \rangle$ is the inner product between two vectors. Then the dual problem associated with problem~\eqref{pb:problem} can be defined as
$\max_{\lambda \in \mathbb{R}^p} \Phi(\lambda)$.
Suppose $x^\star$ is an optimal solution of the primal problem~\eqref{pb:problem} and $\lambda^\star$ is an optimal solution of the dual problem. Then, we make the following assumption:
\begin{assumption}
	\label{assum:strong_dual}
	Strong duality holds for the problem~\eqref{pb:problem}. That is, $f(x^\star) = \Phi(\lambda^\star)$.	
\end{assumption}
Assumption~\ref{assum:strong_dual} implies that $(x^\star, \lambda^\star)$ is a saddle point of the Lagrangian function $L(x;\lambda)$  \citep*{ruszczynski2006nonlinear}. That is, $\forall x \in \mathcal{X}, \lambda \in \mathbb{R}^p$,
\begin{equation}
	\label{eq:saddle_point}
	L(x^\star; \lambda) \leq L(x^\star; \lambda^\star) \leq L(x; \lambda^\star).
\end{equation}
The Augmented Lagrangian function of the primal problem~\eqref{pb:problem} is defined by
$ L_\rho(x; \lambda) = f(x) + \langle Ax-b, \lambda \rangle + \frac{\rho}{2}\| Ax - b \|^2 $,
where $\rho$ is the penalty parameter and $\|\cdot\|$ is the Euclidean norm of a vector. Moreover, we define the augmented dual function
$\Phi_\rho(\lambda) = \min_{x \in \mathcal{X}} L_\rho(x; \lambda)$. Note that $\Phi_\rho(\lambda)$ is always differentiable with respect to $\lambda$ and its gradient is $\nabla \Phi_\rho(\lambda) = A x_\lambda^\star - b$,
where $x_\lambda^\star = \arg\min_{x\in\mathcal{X}} L_\rho(x; \lambda)$. Moreover, $\nabla \Phi_\rho(\lambda)$ is Lipschitz continuous with constant $L_\Phi = \frac{1}{\rho}$ \citep*{rockafellar1976augmented}. The ALM can be viewed as a gradient ascent method on the multiplier $\lambda$ with step size $\frac{1}{L_\Phi} = \rho$. We present the ALM in Algorithm~\ref{alg:exactAL}. Discussion on the convergence of Algorithm~\ref{alg:exactAL} can be found in \cite{rockafellar1976augmented,ruszczynski2006nonlinear} and the references therein. 
ALM converges faster than dual method when the problem is not strongly convex due to its smoothing effect on the dual objective function $\Phi(\lambda)$, \citep{rockafellar1976augmented}.

In practice, Algorithm~\ref{alg:exactAL} cannot be implemented exactly on a fixed-point platform. We present the modified ALM in Algorithm~\ref{alg:inexactAL} to include fixed-point arithmetic errors.
In Algorithm~\ref{alg:inexactAL}, $x^\star_k = \arg \min_{x \in \mathcal{X}} L_\rho(x; \lambda_k)$.
The effect of the fixed-point arithmetic is incorporated in the error terms $\epsilon_{in}^k$ and $\epsilon_{out}^k$. 
Moreover, $K_{out}$ is the number of iterations in the ALM and $L$ is the step size used in the dual update and is defined in Lemma~\ref{lem:inexact_oracle}. Finally, $D$ is a convex and compact set containing at least one optimal multiplier $\lambda^\star$ and $\Pi_D$ denotes the projection onto the set $D$. 
This projection step is the major difference that makes the analysis in this paper different from other works on inexact ALM, e.g., \cite{nedelcu2014computational,necoara2017complexity}.	
The set $D$ is predetermined before running the algorithm. We make the following assumption on $D$:
\begin{assumption}
	\label{assum:choice_D}
	The set $D$ is a box that contains $0$, $2\lambda^\star$ and $\lambda^\star + \mathbf{1}$, where $\mathbf{1}$ is a vector of appropriate dimension and its entries are all $1$.
\end{assumption}

The above choice of $D$ is discussed in more details in Section~\ref{sec:converg}. Note that this set $D$ depends on a uniform bound on $\|\lambda^\star\|$ over a set of problem data. The methods proposed in \citet*{ruszczynski2006nonlinear,mangasarian1985computable,nedic2009approximate,devolder2012double} establish such bounds on $\|\lambda^\star\|$ for fixed problem parameters. On the other hand, \citet*{patrinos2014accelerated,richter2011towards} provide uniform bounds on $\|\lambda^\star\|$ assuming $ b \in \mathcal{B}$ in the constraints in problem~\eqref{pb:problem}. However, the methods in \cite{patrinos2014accelerated,richter2011towards} can only be applied to quadratic problems. For general problems considered in this paper, the interval analysis method in \citet*{hansen1993bounds} can be used to estimate a uniform bound on $\|\lambda^\star\|$. However, the interval analyisis method requires restrictive assumptions on the set of problem data and gives impractical bounds \citep*{jerez2015low}. In practice, we can approximate these bounds using sampling and then apply an appropriate scaling factor. Our contribution in this paper  is to develop a new projected ALM method that employs such bounds to control data overflow on fixed-point platforms.

\section{Convergence Analysis}
\label{sec:converg}
\begin{algorithm}[t]
	\caption{Augmented Lagrangian Method under inexactness}
	\label{alg:inexactAL}
	\begin{algorithmic}[1]
		\Require{Initialize $\lambda_0 \in \mathbb{R}^p$, k=0}
		\While{$k \leq K_{out}$}
		\State{$\tilde{x}_k \approx \arg \min_{x \in \mathcal{X}} L_\rho(x; \lambda_k)$ so that:
			\setlength{\belowdisplayskip}{10pt}
			\setlength{\abovedisplayskip}{10pt}
			\begin{equation*}
				L_\rho(\tilde{x}_k, \lambda_k) - L_\rho(x^\star_k, \lambda_k) \leq \epsilon_{in}^k
			\end{equation*}}
			\State{$\lambda_{k+1} = \Pi_{D}[\lambda_k + \frac{1}{L}(A \tilde{x}_k - b + \epsilon_{out}^k)]$}
			\Let{$k$}{$k + 1$}
			\EndWhile
		\end{algorithmic}
	\end{algorithm}
	
	In this section, we show convergence of Algorithm~\ref{alg:inexactAL} to a neighborhood of the optimal dual and primal objective value. First, we make some necessary assumptions on the boundedness of errors appearing in the algorithm.
	\begin{assumption}
		\label{assum:bound_var}
		At every iteration of Algorithm~\ref{alg:inexactAL}, the errors are uniformly bounded, i.e., 
		\begin{equation*}
			\label{eq:bd_epsilon}
			0 \leq \epsilon_{in}^k \leq B_{in}, \quad \|\epsilon_{out}^k\| \leq B_{out}, \text{ for all } k.
		\end{equation*}
	\end{assumption}
	This assumption is satisfied by selecting appropriate data precision and subproblem solver parameters. Specifically, $\epsilon_{in}^k \geq 0$ if $\tilde{x}_k$ is always feasible, that is, $\tilde{x}_k \in \mathcal{X}$, which is possible if, e.g., the projected gradient method is used to solve the subproblem in line 2 in Algorithm~\ref{alg:inexactAL}. 
	Due to 
	the projection operation in line 3 of Algorithm~\ref{alg:inexactAL}, we have that for $\forall \lambda_1, \lambda_2 \in D_{\delta}$, $\|\lambda_1 - \lambda_2\| \leq B_{\lambda}$, 
	where $D_\delta = D\cup\{\lambda_\delta \in \mathbb{R}^m : \lambda_\delta = \lambda + \epsilon_{out}^k, \text{ for all }\lambda \in D\}$. Since $D$ is compact and $\|\epsilon_{out}^k\|$ is bounded according to Assumption~\ref{assum:bound_var}, $D_\delta$ is also compact. $D_\delta$ is only introduced for the analysis. In practice, we only need to compute the size of $D$ to implement Algorithm~\ref{alg:inexactAL}.

	\subsection{Inexact Oracle}

	Consider the concave function $\Phi_\rho(\lambda)$ with $L_\Phi$-Lipschitz continuous gradient. For any $\lambda_1$ and $\lambda_2 \in \mathbb{R}^p$, we have
	$0 \geq \Phi_\rho(\lambda_2) - [\Phi_\rho(\lambda_1) +  \langle \nabla\Phi_\rho(\lambda_1), \lambda_2 - \lambda_1 \rangle] \geq -\frac{L_\Phi}{2}\|\lambda_1 - \lambda_2\|^2$.
	Recall the expression of $\nabla \Phi_\rho(\lambda)$.
	Since step 2 of Algorithm~\ref{alg:inexactAL} can only be solved approximately, $\nabla \Phi_\rho(\lambda)$ can only be evaluated inexactly.
	Therefore, the above inequalities can not be satisfied exactly.
	We extend the results in \cite{devolder2014first, necoara2017complexity} and propose the following inexact oracle to include the effect of $\epsilon_{out}^k$.

	\begin{lemma}
		\label{lem:inexact_oracle}
		(Inexact Oracle) Let assumptions~\ref{assum:prob_form}, \ref{assum:strong_dual} and \ref{assum:bound_var} hold. Moreover, consider the approximations $\Phi_{\delta,L}(\lambda_k) = L_\rho(\tilde{x}_k;\lambda_k) + B_{out}B_{\lambda} $ to $\Phi_\rho(\lambda_k)$ and $ s_{\delta, L}(\lambda_k) = A \tilde{x}_k - b + \epsilon_{out}^k $
		to $\nabla \Phi_\rho(\lambda_k)$. Then these approximations consititute a $(\delta, L)$ inexact oracle to the concave function $\Phi_\rho(\lambda)$ in the sense that, for $\forall \lambda \in D_{\delta}$,
		\vspace{-1mm}
		\begin{equation}
			\label{eq:inexact_oracle}
			\begin{split}
				0 \geq \Phi_\rho(\lambda) - (\Phi_{\delta, L}(\lambda_k) + & \langle s_{\delta, L}(\lambda_k), \lambda - \lambda_k \rangle) \\
				& \geq -\frac{L}{2}\|\lambda - \lambda_k\|^2 - \delta,
			\end{split}
		\end{equation}
		where $L = 2L_\Phi = \frac{2}{\rho}$ and $\delta = 2B_{in} + 2B_{out}B_\lambda$.
	\end{lemma}
	
	\begin{proof}
		The proof is similar to \cite{necoara2017complexity} and therefore is omitted.
	\end{proof}
	
	\subsection{Dual Suboptimality}
	
	Showing the convergence of the dual variable is similar to showing the convergence of the projected gradient method with the inexact oracle used in \cite{patrinos2015dual}. Therefore, we omit the proof and directly summarize the dual suboptimality results of our algorithm. Specifically, we have the following inequality:
	\begin{equation}
		\label{proof:dual_optimality_5}
		\begin{split}
			\Phi_\rho(\lambda^\star) & - \Phi_\rho(\lambda_{k+1}) \\
			& \leq \frac{L}{2}(\|\lambda_k - \lambda^\star\|^2 - \|\lambda_{k+1} - \lambda^\star\|^2) + \delta. \\
		\end{split}
	\end{equation}
	Furthermore, similar to the Theorem 5 in \cite{patrinos2015dual}, we have the following convergence result for the dual variable:
	\begin{theorem}
		Let assumptions~\ref{assum:prob_form},  \ref{assum:strong_dual} and \ref{assum:bound_var} hold. Define $\bar{\lambda}_K = \frac{1}{K}\sum_{k = 1}^{K}\lambda_k$. Then, we have
		$\Phi_\rho(\lambda^\star) - \Phi_\rho(\bar{\lambda}_k) \leq \frac{L}{2k}\|\lambda_0 - \lambda^\star\|^2 + \delta$.
	\end{theorem}

	\subsection{Primal Infeasibility and Suboptimality}
	
	Define the Lyapunov/Merit function
	$\phi^k(\lambda) = \frac{L}{2}\|\lambda_k - \lambda\|^2 + \frac{1}{2}\|\lambda_{k-1} - \lambda^\star\|^2$.
	Also define the residual function $r(x) = Ax - b$. We have the following intermediate result:
	\begin{lemma}
		\label{lem:intermediate_lemma}
		Let assumptions~\ref{assum:prob_form},\ref{assum:strong_dual} and \ref{assum:bound_var} hold. For all $k \geq 1$, and for all $\lambda \in D$, we have that
		\begin{equation}
			\label{eq:intermediate_lemma}
			f(\tilde{x}_k) - f(x^\star) + \langle \lambda, r(\tilde{x}_k) \rangle \leq \phi^k(\lambda) - \phi^{k+1}(\lambda) + E,
		\end{equation} 
		where $E = (1+\frac{4}{L})B_\lambda B_{out} + (1 + \frac{4}{L})B_{in} + (\frac{1}{2} + \frac{1}{2L})B_{out}^2$.
	\end{lemma}
	
	\begin{proof}
		We recall that $\tilde{x}_k$ is a suboptimal solution as defined in line 2 of Algorithm~\ref{alg:inexactAL} that satisfies
		$ L_\rho(\tilde{x}_k;\lambda_k) - L_\rho(x_k^\star;\lambda_k) \leq \epsilon_{in}^k$.	
		Moreover, due to the optimality of $x_k^\star$, we also have that
		$L_\rho(x_k^\star;\lambda_k) \leq L_\rho(x^\star;\lambda_k) = f(x^\star)$,
		where $x^\star$ is the optimal solution to problem~\eqref{pb:problem}. The equality is because $r(x^\star) = Ax^\star - b = 0$. Combining these two inequalities, we have that
		$L_\rho(\tilde{x}_k;\lambda_k) - f(x^\star) \leq \epsilon_{in}^k$.
		Expanding $L_\rho(\tilde{x}_k;\lambda_k)$ and rearranging terms, we get
		$ f(\tilde{x}_k) - f(x^\star) \leq -\langle \lambda_k, r(\tilde{x}_k) \rangle - \frac{\rho}{2}\|r(\tilde{x}_k)\|^2 + \epsilon_{in}^k$.	
		Adding $\langle \lambda, r(\tilde{x}_k) \rangle$ to both sides of the above inequality, we get
		\setlength{\belowdisplayskip}{10pt}
		\setlength{\abovedisplayskip}{10pt}	
		\begin{equation}
			\label{eq:lemma6_5}
			\begin{split}
				f(\tilde{x}_k) - & f(x^\star) + \langle \lambda, r(\tilde{x}_k) \rangle \\ & \leq \langle \lambda - \lambda_k ,  r(\tilde{x}_k) \rangle - \frac{\rho}{2}\|r(\tilde{x}_k)\|^2 + \epsilon_{in}^k.
			\end{split}		
		\end{equation}	
		In what follows, we show that the right hand side of \eqref{eq:lemma6_5} is upper bounded by $\phi^k(\lambda) - \phi^{k+1}(\lambda) + E$. First, we focus on the term $\langle \lambda - \lambda_k ,  r(\tilde{x}_k) \rangle$. For all $\lambda \in D$, we have that
		$ \|\lambda_{k+1} - \lambda\|^2 = \|\Pi_D[\lambda_k + \frac{1}{L}(r(\tilde{x}_k) + \epsilon_{out}^k)] - \lambda\|^2  \nonumber \leq \|\lambda_k + \frac{1}{L}(r(\tilde{x}_k) + \epsilon_{out}^k)  - \lambda\|^2 = \|\lambda_k - \lambda\|^2 + 2\langle \frac{1}{L}(r(\tilde{x}_k) + \epsilon_{out}^k), \lambda_k - \lambda \rangle + \frac{1}{L^2}\|r(\tilde{x}_k) + \epsilon_{out}^k\|^2 $,	
		where the inequality follows from the contraction of the projection onto a convex set. Rearranging terms in the above inequality and multiplying both sides by $\frac{L}{2}$, we have
		$ \langle r(\tilde{x}_k), \lambda - \lambda_k \rangle \leq \frac{L}{2}(\|\lambda_k - \lambda\|^2 - \|\lambda_{k+1} - \lambda\|^2)
		- \langle \epsilon_{out}^k, \lambda - \lambda_k \rangle
		+ \frac{1}{2L}\|r(\tilde{x}_k)\|^2 + \frac{1}{L}\langle r(\tilde{x}_k), \epsilon_{out}^k \rangle + \frac{1}{2L}\|\epsilon_{out}^k\|^2$.
		Applying Assumption~\ref{assum:bound_var} and the Cauchy-Schwartz inequality, we obtain
		\begin{flalign}
			\label{eq:lemma6_7}
			& \langle r(\tilde{x}_k), \lambda - \lambda_k \rangle \leq
			\frac{L}{2}(\|\lambda_k - \lambda\|^2 - \|\lambda_{k+1} - \lambda\|^2) & \\
			& +B_{out}B_\lambda + \frac{1}{2L}\|r(\tilde{x}_k)\|^2 + \frac{1}{L}\langle r(\tilde{x}_k), \epsilon_{out}^k \rangle + \frac{1}{2L}B_{out}^2. & \nonumber
		\end{flalign}	
		To upper bound the term $\langle r(\tilde{x}_k), \epsilon_{out}^k \rangle$ in \eqref{eq:lemma6_7}, first we add and subtract $\epsilon_{out}^k$ from $r(\tilde{x}_k)$, and rearrange terms to get
		$ \langle r(\tilde{x}_k), \epsilon_{out}^k \rangle = \langle r(\tilde{x}_k) + \epsilon_{out}^k, \epsilon_{out}^k \rangle - \|\epsilon_{out}^k\|^2$.
		Defining $\lambda_{\delta k} = \lambda_k + \epsilon_{out}^k$, recalling the definition of $s_{\delta,L}(\lambda_k)$ in Lemma~\ref{lem:inexact_oracle}, and ignoring the term $-\|\epsilon_{out}^k\|^2$, we obtain
		$ \langle r(\tilde{x}_k), \epsilon_{out}^k \rangle \leq \langle s_{\delta,L}(\lambda_k), \lambda_{\delta k} - \lambda_k \rangle$.	
		Next we apply the second inequality in Lemma~\ref{lem:inexact_oracle} to upper bound $\langle s_{\delta,L}(\lambda_k), \lambda_{\delta k} - \lambda_k \rangle$. In order to apply Lemma~\ref{lem:inexact_oracle}, both $\lambda_{\delta k}$ and $\lambda_k$ need to belong to $D_\delta$. Due to the projection in line 3 in Algorithm~\ref{alg:inexactAL}, $\lambda_k$ always belongs to $D$. Recalling the definition of $D_\delta$, it is straightforward to verify that $\lambda_k$, $\lambda_{\delta k} \in D_\delta$. Thus applying Lemma~\ref{lem:inexact_oracle} we get
		\setlength{\belowdisplayskip}{10pt}
		\setlength{\abovedisplayskip}{10pt}
		\begin{equation}
			\label{eq:lemma6_8}
			\begin{split}
				& \langle r(\tilde{x}_k), \epsilon_{out}^k \rangle  \leq \langle s_{\delta,L}(\lambda_k), \lambda_{\delta k} - \lambda_k \rangle \\
				& \leq \Phi_\rho(\lambda_{\delta k}) - \Phi_{\delta,L}(\lambda_k) + \frac{L}{2}B_{out}^2 + \delta.
			\end{split}
		\end{equation}
		Since $\lambda^\star$ is the global maximizer of the function $\Phi_\rho(\lambda)$, we get $\Phi_\rho(\lambda^\star) \geq \Phi_\rho(\lambda_{\delta k})$. We can also show that $\Phi_{\delta, L}(\lambda_k) \geq \Phi_\rho(\lambda_k)$ always holds because
		$ \Phi_{\delta,L}(\lambda_k) = L_\rho(\tilde{x}_k;\lambda_k) + B_{out}B_\lambda \geq L_\rho(\tilde{x}_k;\lambda_k) \geq L_\rho(x_k^\star;\lambda_k)$. Combining these two inequalities we obtain that $\Phi_\rho(\lambda^\star) - \Phi_\rho(\lambda_k) \geq \Phi_\rho(\lambda_{\delta k}) - \Phi_{\delta,L}(\lambda_k)$. Substituting this inequality into \eqref{eq:lemma6_8}, we have that
		$ \langle r(\tilde{x}_k), \epsilon_{out}^k \rangle \leq \Phi_\rho(\lambda^\star) - \Phi_\rho(\lambda_k) + \frac{L}{2}B_{out}^2 + \delta $.
		Combining this inequality, \eqref{proof:dual_optimality_5} and \eqref{eq:lemma6_7}, we get
		\begin{equation}
			\label{eq:lemma6_11}
			\begin{split}
				& \langle r(\tilde{x}_k), \lambda - \lambda_k \rangle \leq \frac{L}{2}(\|\lambda_k - \lambda\|^2 - \|\lambda_{k+1} - \lambda\|^2) \\
				& + \frac{1}{2}(\|\lambda_{k-1} - \lambda^\star\|^2 - \|\lambda_{k} - \lambda^\star\|^2) + \frac{1}{2L}\|r(\tilde{x}_k)\|^2 \\
				& + B_{out}B_\lambda + (\frac{1}{2} + \frac{1}{2L})B_{out}^2 + \frac{2}{L}\delta.
			\end{split}
		\end{equation}
		Combining \eqref{eq:lemma6_11} with \eqref{eq:lemma6_5}, we have that
		$ f(\tilde{x}_k) - f(x^\star) + \langle \lambda, r(\tilde{x}_k) \rangle \leq \frac{L}{2}(\|\lambda_k - \lambda\|^2 - \|\lambda_{k+1} - \lambda\|^2) + \frac{1}{2}(\|\lambda_{k-1} - \lambda^\star\|^2 - \|\lambda_{k} - \lambda^\star\|^2)  + (\frac{1}{2L} - \frac{\rho}{2})\|r(\tilde{x}_k)\|^2 + B_{out}B_\lambda + (\frac{1}{2} + \frac{1}{2L})B_{out}^2 + \frac{2}{L}\delta + \epsilon_{in}^k$.
		From Lemma~\ref{lem:inexact_oracle}, we have that $L = 2L_\Phi = \frac{2}{\rho}$, and $\frac{1}{2L} - \frac{\rho}{2} < 0$. Therefore, the term $(\frac{1}{2L} - \frac{\rho}{2})\|r(\tilde{x}_k)\|^2 < 0$ can be neglected. Recalling the definition of $\phi^k(\lambda)$, $\delta$ and $E$, we obtain \eqref{eq:intermediate_lemma}, which completes the proof.
	\end{proof}
	Next, we apply Lemma~\ref{lem:intermediate_lemma} to prove the primal suboptimality and infeasibility of Algorithm \ref{alg:inexactAL}.
	%
	\begin{theorem}
		\label{thm:primal_opt_feas}
		(Primal Suboptimality and Infeasibility) Let assumptions~\ref{assum:prob_form}, \ref{assum:strong_dual} and \ref{assum:bound_var} hold. Define $\bar{x}_K = \frac{1}{K}\sum_{k=1}^{K}\tilde{x}_k$. Then, we have that (a) primal optimality:
		\begin{equation}
			\label{eq:primal_opt_1}
			-(\frac{1}{K}\phi^1(2\lambda^\star) + E) \leq f(\bar{x}_K) - f(x^\star) \leq \frac{1}{K}\phi^1(0) + E,
		\end{equation}
		\noindent (b) primal feasibility:
		\begin{equation}
			\label{eq:primal_feas}
			\|r(\bar{x}_K)\| \leq \frac{1}{K}\phi^1(\lambda^\star + \frac{r(\bar{x}_K)}{\|r(\bar{x}_K)\|}) + E.
		\end{equation}
	\end{theorem}
	
	\begin{proof}
		Summing inequality (\ref{eq:intermediate_lemma}) in Lemma~\ref{lem:intermediate_lemma} for $k = 1, 2, \dots, K$, we have that
		$ \sum_{k=1}^{K}  f(\tilde{x}_k) - Kf(x^\star) + \langle \lambda, \sum_{k=1}^{K}r(\tilde{x}_k) \rangle \leq \phi^1(\lambda) - \phi^{K+1}(\lambda) + KE \leq \phi^1(\lambda) + KE$, 	
		where the second inequality follows from $\phi^k(\lambda) \geq 0$. Dividing both sides of the above inequality by $K$ and using the fact that $\frac{1}{K}\sum_{k=1}^{K}r(\tilde{x}_k) = r(\bar{x}_K)$, we get
		$\label{proof:final_2}
		\frac{1}{K}\sum_{k=1}^{K}f(\tilde{x}_k) - f(x^\star) + \langle \lambda, r(\bar{x}_K) \rangle \leq \frac{1}{K}\phi^1(\lambda) + E$.
		From the convexity of the function $f(x)$, we have that $f(\bar{x}_K) \leq \frac{1}{K}\sum_{k=1}^{K}f(\tilde{x}_k)$. Combining the above two inequalities, we obtain
		\begin{equation}
			\label{proof:final_3}
			f(\bar{x}_K) - f(x^\star) + \langle \lambda, r(\bar{x}_K) \rangle \leq \frac{1}{K}\phi^1(\lambda) + E.	
		\end{equation}	
		To prove the right inequality in \eqref{eq:primal_opt_1}, let $\lambda = 0$ in \eqref{proof:final_3}. Then, we have that
		$ f(\bar{x}_K) - f(x^\star) \leq \frac{1}{K}\phi^1(0) + E$.		
		To prove the left inequality in \eqref{eq:primal_opt_1}, according to the inequality~\eqref{eq:saddle_point} and Assumption~\ref{assum:strong_dual}, we have that
		\begin{equation}
			\label{proof:final_4_1}
			f(x^\star) - f(\bar{x}_K) \leq \langle \lambda^\star, r(\bar{x}_K) \rangle.
		\end{equation}	
		Adding $\langle \lambda^\star, r(\bar{x}_K) \rangle$ to both sides of \eqref{proof:final_4_1} and rearranging terms, we have that
		\begin{equation}
			\label{proof:final_4_2}
			\langle \lambda^\star, r(\bar{x}_K) \rangle \leq f(\bar{x}_K) - f(x^\star) + \langle 2\lambda^\star, r(\bar{x}_K) \rangle.
		\end{equation}
		Combining inequalities (\ref{proof:final_4_1}) and (\ref{proof:final_4_2}), we obtain
		$ f(x^\star) - f(\bar{x}_K) \leq f(\bar{x}_K) - f(x^\star) + \langle 2\lambda^\star, r(\bar{x}_K) \rangle$.
		Combining this inequality and the inequality in~\eqref{proof:final_3} with $\lambda = 2\lambda^\star$, we get
		$ f(x^\star) - f(\bar{x}_K) \leq \frac{1}{K}\phi^1(2\lambda^\star) + E$.
		
		To prove \eqref{eq:primal_feas}, let $\lambda = \lambda^\star + \frac{r(\bar{x}_K)}{\|r(\bar{x}_K)\|}$ in \eqref{proof:final_3}, which also belongs in $D$ by Assumption~\ref{assum:choice_D}. Rearranging terms, we obtain
		$f(\bar{x}_K) - f(x^\star) +  \langle \lambda^\star, r(\bar{x}_K) \rangle + \|r(\bar{x}_K)\| \leq \frac{1}{K}\phi^1(\lambda^\star + \frac{r(\bar{x}_K)}{\|r(\bar{x}_K)\|}) + E$. 
		Letting $x = \bar{x}_K$ in the second inequality in \eqref{eq:saddle_point}, we have that $f(\bar{x}_K) + \langle \lambda^\star, r(\bar{x}_K) \rangle - f(x^\star) \geq 0$. Combining this inequality with the above inequality, we get
		$\|r(\bar{x}_K)\| \leq \frac{1}{K}\phi^1(\lambda^\star + \frac{r(\bar{x}_K)}{\|r(\bar{x}_K)\|}) + E$,
		which completes the proof.
	\end{proof}

\section{Fixed-point Implementation}
\label{sec:fp}
To apply the Algorithm~\ref{alg:inexactAL} to solve the problem~\eqref{pb:problem} on fixed-point platforms, we need another assumption.

\begin{assumption}
	\label{assum:4.1}
	The set $\mathcal{X}$ is compact and simple. 
\end{assumption}

The compactness is needed to bound the primal iterates $\tilde{x}_k$ in Algorithm~\ref{alg:inexactAL}.
$\mathcal{X}$ is simple in the sense that computing the projection onto $\mathcal{X}$ is of low complexity on fixed-point platforms, for example, $\mathcal{X}$ is a box. 
We note that our convergence analysis in Section~\ref{sec:converg} does not rely on the Assumption~\ref{assum:4.1}. This assumption is needed only for fixed-point implementation.
Due to the projection operation in the multiplier update in Algorithm~\ref{alg:inexactAL}, $\|\lambda_k\|$ is always bounded. Since now both $x_k$ and $\lambda_k$ lie in compact sets, the bounds on the norms of these variables can be used to design the word length of the fixed-point data to avoid overflow, similar to Section 5.3 in \cite{patrinos2015dual}. Meanwhile, according to the results in \eqref{eq:primal_opt_1} and \eqref{eq:primal_feas}, we can choose the number $K_{out}$ of the outer iterations, the accuracy of the multiplier update $B_{out}$ and the accuracy of the subproblem solution $B_{in}$ in Algorithm~\ref{alg:inexactAL} to achieve an $\epsilon$-solution to the problem~\eqref{pb:problem}. 

Specifically, to achieve the accuracy $B_{in}$, we can determine the iteration complexity of the subproblem solver based on the studies in \cite{devolder2014first,patrinos2015dual, schmidt2011convergence}. However, these results are usually conservative in practice.
Therefore, in this section, we propose a stopping criterion for early termination of the subproblem solver with guarantees on the solution accuracy. This stopping criterion can be used for box-constrained convex optimization problems that are not necessarily strongly convex, and is simple  to check on embedded devices.
A similar stopping criterion has been proposed in \cite{nedelcu2014computational}, which however requires the objective function of the subproblem to be strongly convex. In what follows, we relax the strong convexity assumption so that this stopping criterion can be applied.
\begin{lemma}
	\label{lem:stopcondition}
	Consider a convex and nonsmooth function $f(x)$, with support $\mathcal{X} \subset \mathbb{R}^n$ and set of minimizers $X^\star$. Assume $f(x)$ satisfies the quadratic growth condition, 
	\begin{equation}
		\label{eq:qgc_1}
		f(x) \geq f^\star + \frac{\sigma}{2}\mathtt{dist}^2(x, X^\star), \text{ for } \forall x \in \mathcal{X},
	\end{equation}
	where $\mathtt{dist}(x, X^\star)$ is the distance from $x$ to the set $X^\star$. and $\sigma>0$ is a scalar. Then, we have that for any $x \in \mathcal{X}$ and any $s \in \mathcal{N}_{\mathcal{X}}(x)$,
	\begin{equation}
		\label{eq:qgc_2}
		f(x) - f^\star \leq \frac{2}{\sigma}\|\partial f(x) + s\|^2,
	\end{equation}
	where $\mathcal{N}_\mathcal{X}(x)$ is the normal cone at $x$ with respect to $\mathcal{X}$.
\end{lemma}

\begin{proof}
	The proof is the similar to the proof in \cite{nedelcu2014computational} and therefore is omitted.
\end{proof}

\citet*{necoara2018linear} have shown that the function $h(Ax)$ satisfies the inequality~\eqref{eq:qgc_1} when the set $\mathcal{X}$ is polyhedral and $h(y)$ is strongly convex. However, since the subproblem objective $L_\rho(x,\lambda_k)$ is in the form $\sum_i h_i(A_ix)$ rather than $h(Ax)$, in what follows we present a generalization of Theorem 8 in \cite{necoara2018linear} so that the function $L_\rho(x,\lambda_k)$ satisfies the inequality~\eqref{eq:qgc_1} and the condition~\eqref{eq:qgc_2} can be applied.

\begin{lemma}
	\label{lem:qg_func}
	Consider a convex function $H(x) = \sum_i h_i(A_ix - b_i)$ with the polyhedral support $\mathcal{X} = \{ x:Cx \leq d \}$, where $h_i(y)$ is $\sigma_i-$strongly convex with respect to $y$ for any $i$. Then we have
	\begin{equation}
		\label{eq:qg_func}
		H(x) \geq H^\star + \frac{\sigma}{2}\mathtt{dist}^2(x, X^\star),
	\end{equation}
	where $\sigma = \frac{\min_i\{\sigma_i\}}{\theta^2(\tilde{A},C)}$, the matrix $\tilde{A}$ is in the form
	$[\cdots | A_i^T | \cdots ]^T$
	and $\theta(\tilde{A}, C)$ is a constant only related to the matrices $\{A_i\}$ and $C$.
\end{lemma}
\begin{proof}
	Given $x^\star \in X^\star$, for any $x \in \mathcal{X}$, from the strong convexity of the function $h_i$, we have that
	$ h_i(A_ix - b_i) \geq  h_i(A_ix^\star-b_i) + \langle \partial h_i(y)|_{y = A_ix^\star - b_i}, A_i(x - x^\star) \rangle + \frac{\sigma_i}{2}\|A_ix - A_ix^\star\|^2 = h_i(A_ix^\star-b_i) + \langle  x - x^\star, A_i^T \partial h_i(y)|_{y = A_ix^\star - b_i} \rangle + \frac{\sigma_i}{2}(x - x^\star)^TA_i^TA_i(x - x^\star)$.
	Adding up the above inequalities for all $i$, we have
	$H(x) \geq  H(x^\star) + \langle \sum_i A_i^T\partial h_i(y)|_{y = A_ix^\star - b_i}, x - x^\star \rangle  + \sum_i \frac{\sigma_i}{2}(x - x^\star)^TA_i^TA_i(x - x^\star)$.
	Recall that $\partial H(x^\star) = \sum_i A_i^T\partial h_i(y)|_{y = A_ix^\star - b_i}$. From the optimality of $x^\star$, we have that $\langle \sum_i A_i^T\partial h_i(y)|_{y = A_ix^\star - b_i}, x - x^\star \rangle \geq 0$ for any $x \in \mathcal{X}$. Therefore, we obtain
	\begin{equation}
		\label{eq:qgc_func_1}
		H(x) \geq H(x^\star) + \sum_i \frac{\sigma_i}{2}(x - x^\star)^TA_i^TA_i(x - x^\star).
	\end{equation}
	Next, we show that the optimal solution set $X^\star = \{\tilde{x}^\star: \tilde{A}\tilde{x}^\star = t^\star ; C\tilde{x}^\star \leq d \}$, where $t^\star = \tilde{A}x^\star$. To do so, we decompose the space $\mathbb{R}^n$ into two mutually orthogonal subspaces, $\mathcal{V}$ and $\mathcal{U}$, where $\mathcal{V}$ is the intersection of the kernel spaces of all matrices $A_i$ and $\mathcal{U}$ is the union of the row space of all matrices $A_i$. Showing that ${X}^\star = \{\tilde{x}^\star: \tilde{A}\tilde{x}^\star = t^\star ; C\tilde{x}^\star \leq d \}$ is equivalent as showing that $X^\star = \tilde{X}^\star$, where $\tilde{X}^\star = \{\tilde{x}^\star: \tilde{x}^\star = x^\star + v, \text{ where } v \in \mathcal{V} ; C\tilde{x}^\star \leq d \}$. First, we show that $\tilde{X}^\star \subset X^\star$. If $\tilde{x}^\star \in \tilde{X}^\star$, we have that $\tilde{x}^\star$ is feasible and $H(\tilde{x}^\star) = H(x^\star)$ because  $A_i\tilde{x}^\star = A_i(x^\star + v) = A_ix^\star$, for any $i$. Therefore, $\tilde{X}^\star \subset X^\star$. Second, we prove that $X^\star \subset \tilde{X}^\star$ by showing that if $\tilde{x} \notin \tilde{X}^\star$, then $\tilde{x}$ is not optimal. If $C\tilde{x} > d$, then $\tilde{x}$ is not feasible. On the other hand, if $C\tilde{x} \leq d$ but $\tilde{x} = x^\star + \alpha v + u$, where $\alpha$ is any real number and $u \in \mathcal{U}$, then $H(\tilde{x}) > H(x^\star)$. This is due to the inequality~\eqref{eq:qgc_func_1} and $\sum_i \frac{\sigma_i}{2}(\tilde{x} - x^\star)^TA_i^TA_i(\tilde{x} - x^\star) > 0$. Therefore, $X^\star \subset \tilde{X}^\star$ and we have shown that $X^\star = \tilde{X}^\star$.	
	Since $X^\star = \{\tilde{x}^\star: \tilde{A}\tilde{x}^\star = t^\star ; C\tilde{x}^\star \leq d \}$, according to Lemma 15 in \cite{wang2014iteration}, for any $x \in \mathcal{X}$, we have that $\mathtt{dist}(x, X^\star) \leq \theta(\tilde{A},C) \|\tilde{A}(x - x^\star)\|$. Furthermore, we have that $\mathtt{dist}^2(x, X^\star) \leq \theta^2(\tilde{A},C)\|\tilde{A}(x - x^\star)\|^2 \leq  \frac{\theta^2(\tilde{A},C)}{\min_i\{\sigma_i\}}\sum_i \sigma_i(x - x^\star)^TA_i^TA_i(x - x^\star)$. Combining this inequality with \eqref{eq:qgc_func_1}, we get the desired result \eqref{eq:qg_func}.	
\end{proof}

The constant $\theta(\tilde{A},C)$ can be obtained as in \citet*{necoara2018linear}.
If the objective function $f(x)$ in problem~\eqref{pb:problem} is of the form $h(Ax - b)$ where $h(y)$ is strongly convex and the set  $\mathcal{X} = \{x:l \leq x \leq u\}$, according to Lemma  \ref{lem:qg_func}, the subproblem objective $L_\rho(x, \lambda_k)$ satisfies the inequality~\eqref{eq:qgc_1}. Therefore,  according to Lemma~\ref{lem:stopcondition}, we have that the iterate $x_t$ returned from the subproblem solver satisfies $L_\rho(x_t;\lambda_k) - L_\rho(x_k^\star;\lambda_k) \leq B_{in}$ if $ \|\partial L_\rho(x_t;\lambda_k) + s_t^\star\| \leq \sqrt{\frac{\sigma}{2}B_{in}}$, where $s_t^\star = \arg\min_{s_t \in \mathcal{N}_{\mathcal{X}}(x_t)}\|\partial L_\rho(x_t;\lambda_k) + s_t\|$. As discussed in \cite{nedelcu2014computational}, $\|\partial L_\rho(x_t;\lambda_k) + s_t^\star\|$ can be efficiently evaluated on the embedded platform.

\section{Numerical Simulations}
\label{sec:sim}
\vspace{-1mm}
In this section, we present simulation results for a utility maximization example to verify the convergence and error analysis results in Section~\ref{sec:converg} and the design of the fixed-point implementation in Section~\ref{sec:fp}. The simulations are conducted using the Fixed-Point Designer in Matlab R2015a on a Macbook Pro with 2.6GHz Intel Core i5 and 8GB, 1600MHz memory.
Consider an undirected graph $G = (\mathcal{N}, \mathcal{E})$, where $\mathcal{N} = \{1,2,\dots,N\}$ is the set of nodes and $\mathcal{E}$ is the set of edges, so that $(i,j) \in \mathcal{E}$ if the nodes $i$ and $j$ are connected in the graph $G$. Denote the set of neighbors of node $i$ as $\mathcal{N}_i$. The set $\mathcal{N}$ consists of two subsets $\{S, D\}$, where $S$ and $D$ are the sets of source and the destination nodes, respectively. The node $i \in S$ generates data at a rate $s_i$, where $s_{\min} \leq s_i \leq s_{\max}$. The data flows from node $i$ to node $j$ through edge $(i,j) \in \mathcal{E}$ at a rate $t_{ij}$, where $0 \leq t_{ij} \leq c_{ij}$. All generated data finally flows into the destination nodes, which are modeled as sinks and can absorb the incoming data at any rates. The nodes collaboratively solve the following network utility maximization (NUM) problem
\begin{flalign}
	\label{eq:pb_num}
	& \max_{\{s_i\}, \{t_{ij}\}} \; \sum_{i = 1} \log(s_i) & \nonumber \\
	& \text{s.t. } \; \sum_{j \in \mathcal{N}_i} t_{ij} - \sum_{j \in \mathcal{N}_i} t_{ji} = s_i, \;\;\; \forall i \in S & \\
	& s_{\min} \leq s_i \leq s_{\max}, \; \forall i \in S, \;\;\; 0 \leq t_{ij} \leq c_{ij}, \;\;\; \forall (i,j) \in \mathcal{E}, \nonumber
\end{flalign}
where the constraint $\sum_{j \in \mathcal{N}_i} t_{ij} - \sum_{j \in \mathcal{N}_i} t_{ji} = s_i$ expresses the flow conservation law at the node $i$. 
The logarithm objective function is used to measure the utility of the data generation rate.
To solve problem~\eqref{eq:pb_num} distributedly, distributed ALM schemes \citep*{chang2015multi,chatzipanagiotis2015augmented,chatzipanagiotis2016distributed,chatzipanagiotis2017convergence,lee2017complexity} have been proposed that converge much faster than the dual decomposition method although at the cost of solving nontrivial subproblems locally at each iteration. Here we employ the consensus-ADMM method in \cite{chang2015multi} to solve problem~\eqref{eq:pb_num}. Specifically, let node $i$ keep a local decision variable $[s_i, t_{(i)}^T]^T$, where $t_{(i)} = [t_{i1}^{(i)}, \dots, t_{i|\mathcal{N}_i|}^{(i)}, t_{1i}^{(i)}, \dots, t_{|\mathcal{N}_i|i}^{(i)}]^T$. Then, at the $t$ th iteration, each node needs to solve a local problem
\begin{flalign}
	\label{eq:pb_resource_local}
	& \min_{s_i, t_{ij}^{(i)}, t_{ji}^{(i)}} \; -\log(s_i) + \langle p_i^t, t_{(i)} \rangle +  \mu\|t_{(i)} - g_i^t\|^2 & \nonumber \\
	& \text{s.t. } \; \sum_{j \in \mathcal{N}_i} t_{ij}^{(i)} - \sum_{j \in \mathcal{N}_i} t_{ji}^{(i)} = s_i, & \\
	& s_{\min} \leq s_i \leq s_{\max}, \;\;\; 0 \leq t_{ij}^{(i)} \leq c_{ij}, \;\;\; 0 \leq t_{ji}^{(i)} \leq c_{ji}, & \nonumber
\end{flalign} 
where the variables $p_i^t$ and $g_i^t$ are updated at every iteration of the consensus-ADMM to finally achieve consensus $t_{ij}^{(i)} = t_{ij}^{(j)}$ on all edges. Problem \eqref{eq:pb_resource_local} has the form of problem \eqref{pb:problem} and we can apply Algorithm \ref{alg:inexactAL} to solve it using fixed-point data. Since the objective function in \eqref{eq:pb_resource_local} is not quadratic, the results in \cite{jerez2014embedded} cannot be applied.

For our numerical simulations, we randomly generate problem \eqref{eq:pb_num} on a network of $10$ nodes and apply Algorithm~\ref{alg:inexactAL} to solve the subproblem \eqref{eq:pb_resource_local}. 
In what follows we present results for the node that solves the largest subproblems that are of dimension $9$.
The parameters $p_i^t$ and $g_i^t$ in \eqref{eq:pb_resource_local} are obtained by running $30$ iterations of consensus-ADMM on problem~\eqref{eq:pb_num} with double floating point data. This creates $30$ instances of subproblem~\eqref{eq:pb_resource_local} which we solve using Algorithm~\ref{alg:inexactAL}. The bounds on $\|\lambda^\star\|$ for each problem instance were obtained by solving each problem using the Matlab function $fmincon$. Table~\ref{table:primal_feas_opt_num} shows the achieved optimality and feasibility for the worst-case scenario and for three solution accuracies, $\epsilon = 1, 0.1$ and $0.01$. We observe that the theoretical bounds are around 20 to 30 times higher than the actual algorithm performance.
\begin{table}[t]
	\centering
	\scriptsize
	\caption{Primal optimality and feasiblity achieved by Algorithm~\ref{alg:inexactAL} for the NUM problem }
	\label{table:primal_feas_opt_num}	
	\vspace{-4mm} 	
	\begin{tabular}{| c | c | c | c |}
		\multicolumn{4}{c}{} \\
		\hline
		$\epsilon$ & $fl$-$wl$ & low. \eqref{eq:primal_opt_1}, $\;$ opt., $\;$ up. \eqref{eq:primal_opt_1} &  feas., $\;$ \eqref{eq:primal_feas} \\ 
		\hline
		1 & 10-14 & -0.9861, $\;$ 0.0338, $\;$ 0.7026 & 0.0439, $\;$ 1.0 \\
		0.1 & 14-18 & -0.0995, $\;$ 0.0034, $\;$ 0.0707 & 0.0044, $\;$ 0.1 \\
		0.01 & 17-21 & -0.0100, $\;$ 0.00035, $\;$ 0.0071 & 0.00044, $\;$ 0.01 \\	
		\hline
	\end{tabular}	
\end{table}

Note that this simulation serves only the purpose of showing the ability of Algorithm~\ref{alg:inexactAL} to solve non-quadratic convex optimization problems, as the subproblems \eqref{eq:pb_resource_local}. 
The method in Jerez et al. (2014) can also be extended to solve such problems if it is used as the QP solver in the Sequential Quadratic Programming (SQP) framework. However, a new KKT matrix needs to be inverted at each iteration of SQP according to Jerez et al. (2014) and the solutions to the QP problems will be inexact on a fixed-point platform. The effect of this inexactness on the iteration complexity and final solution accuracy of SQP under fixed-point arithmetics has not been theoretically studied.
A systematic solution of \eqref{eq:pb_num} using distributed ALM methods with fixed-point data is an open problem and is left for future research.

\section{Conclusion}
\label{sec:conclude}
In this paper we proposed an Augmented Lagrangian Method to solve convex and non-smooth optimization problems using fixed-point arithmetic. To avoid data overflow, we introduced a projection operation in the multiplier update. Moreover, we present a stopping criterion to terminate early the primal subproblem iteration, while ensuring a desired accuracy of the solution. We presented convergence rate results as well as bounds on the optimality and feasibility gaps.  To the best of our knowledge, this is the first fixed-point ALM that can handle non-smooth problems, data overflow, and can efficiently and systematically utilize iterative solvers in the primal update.

\ifCLASSOPTIONcaptionsoff
  \newpage
\fi

\bibliographystyle{IEEEtranN}
\bibliography{bib_Yan}

%

\end{document}